\documentclass[10pt]{amsart}
\usepackage{amssymb,latexsym,amsmath}
\numberwithin{equation}{section}
\newtheorem{thm}[equation]{Theorem}
\newtheorem{prop}[equation]{Proposition}
\newtheorem{cor}[equation]{Corollary}

\newtheorem{dfn}[equation]{Definition}
\def\t{\tilde}
\def\fS{\mathfrak{S}}
\def\cS{\mathcal{S}}
\gdef\P{\mathbb{P}}
\def\C{\mathbb{C}}

\def\F{{\mathcal{F}}}

\gdef\GSL{{\mathrm{SL}(W)}}

\gdef\gsl{{\mathfrak{sl}(W)}}

\def\l{\ell}
\def\a{{\alpha}}
\def\b{\beta}
\def\c{\gamma}
\def\d{\delta}
\def\e{\varepsilon}
\def\r{\rho}
\def\z{\zeta}

\def\w{\omega}

\def\oa{{\overline \alpha}}
\def\ob{{\overline \beta}}
\def\oc{{\overline \gamma}}

\def\oe{{\overline \varepsilon}}
\def\ovr{{\overline \rho}}
\def\oz{{\overline \zeta}}

\def\half{\tfrac{1}{2}}
\def\fourth{\tfrac{1}{4}}

\def\threehalves{\tfrac{3}{2}}
\begin{document}
\title{The adjoint variety of $\mathrm{SL}_{m+1}\C$ is rigid to order three}
\author{Colleen Robles}
\date{last updated 18 August 2006}
\begin{abstract}
I prove that the adjoint variety of $\mathrm{SL}_{m+1}\C$ in 
$\P(\mathfrak{sl}_{m+1}\C)$ is rigid to order three.
\end{abstract}
\maketitle


The principle result of this paper is Theorem \ref{thm:SL} 
(page \pageref{thm:SL})
which asserts that the adjoint variety of the simple Lie group 
$\mathrm{SL}_{m+1} \C$ is rigid to order three.  The result is extrinsic; 
roughly speaking, if a variety 
$Y \subset \P(\mathfrak{sl}_{m+1}\C) = \P^{m^2 + 2m - 1}$, of dimension 
$n=2m-1$, resembles the adjoint variety to third order at a $3$-general point 
$y\in Y$, then there is a transformation in $\textrm{GL}_{m^2+2m}\C$ 
mapping $Y$ onto the adjoint variety.

The conclusion is significant because it is the first rigidity result for 
a variety with non-vanishing Fubini cubic $F_3$ (a third order invariant).  
And it is striking that this is the first example of $k$-th order rigidity for 
which the $(k+1)$-th order Fubini invariant is nonzero: $F_4$ can not be 
normalized to zero.

The proof is based on the E. Cartan's method of moving frames.  The reader may 
find similar applications of the technique to the study of submanifolds of 
$\C\P^N$ in \cite{GH79, L96, L99, L06, LM02}, and their references.
The paper is organized as follows
\begin{itemize}
  \item[\S\ref{sec:intro}]
    Notation is set.  The first-order adapted frame bundle associated to a 
    variety is introduced, and the relative differential invariants $F_k$, or 
    Fubini forms, are discussed.  (These invariants describe the lines 
    osculating to a variety with order $k$.)
  \item[\S\ref{sec:rigidity}]
    The notions of {\it agreement to $k$-th order} and {\it rigidity} are made 
    explicit, and previous rigidity results are reviewed.
  \item[\S\ref{sec:normalizations}]
    The frame bundle normalizations necessary to our computations are sketched.
  \item[\S\ref{sec:SL}]
    The adjoint variety of $\textrm{SL}_{m+1}\C$ is introduced and the main 
    result, Theorem \ref{thm:SL}, stated.
  \item[\S\ref{sec:proof}]  
    Theorem \ref{thm:SL} is proven.
  \item[\S\ref{sec:conclusions}] An open question: what rigidity results 
    might we expect for the other adjoint varieties?
\end{itemize}
\vspace{0.1in}

\noindent
{\it Acknowledgments.}  I thank J.M. Landsberg for posing the question, and 
several illuminating discussions.  The main result of 
this paper was first presented at the IMA's 2006 summer program on Symmetries 
and Overdetermined Systems of Partial Differential Equations.  I am grateful 
to the AWM for the Travel Grant supporting my attendance. 


\section{Preliminaries}
\label{sec:intro}


\subsection{Notation}
\label{sec:notation}
The aim of this section is merely to establish notation; see \cite[\S3]{IL03}
for a through discussion of the ideas presented here.
I will, for the most part, follow the notation and conventions of 
\cite{IL03, L96}.  

Let $V$ be a complex vector space of dimension $N+1$, and $X^n \subset \P V$ a 
variety of (complex) dimension $n$.  Fix the index ranges
\begin{eqnarray*}
  & 0 \le I, J, K, \ldots \le N & \\
  & 1 \le a , b , c, \ldots \le n & \\
  & n + 1 \le u , v , w \ldots \le N &
\end{eqnarray*}
The Einstein convention holds here: repeated indices, one raised and one 
lowered, are summed over.  Given a collection of vectors $\{ v_j \}$ in a 
vector space, let 
\begin{displaymath}
  \langle v_j \rangle = \hbox{span}\{v_j\}
\end{displaymath}
denote their span.  Given a nonzero $v \in V$, I will denote by $[v]$ its 
projection to $\P V$.  The {\it cone over $X \subset \P V$} is 
\begin{displaymath}
  \widehat X = \{ v \in V \backslash \{0\} \, : \, [v] \in X \} \, .
\end{displaymath}


\subsection{A frame bundle over \boldmath $X$ \unboldmath }

Let $X^n \subset \P V$ be a variety.  The bundle of {\it first-order adapted 
frames} $\F^1$ over 
\begin{displaymath}
M := X_\textit{smooth} \subset X
\end{displaymath}
 is the subset of those $\textbf{e} = (e_0 , e_1 , \ldots , e_N ) \in GL(V)$ 
for which $e_0 \in \widehat M$, and 
the {\it affine tangent space} to the cone over $M$, 
$T_{e_0} \widehat M \subset V$, is spanned by 
$\{ e_0 , e_1 , \ldots , e_n \}$.  Let $\w$ denote the pull-back of the 
Maurer-Cartan form on $GL(V)$ to $\F^1$.  In particular, 
\begin{displaymath}
  d e_I = \w_I^J \, e_J \quad \hbox{ and } \quad
  d \w^I_K = - \w^I_J \wedge \w^J_K \, .
\end{displaymath}
Since the $T_{e_0} \widehat M = \langle e_0 , e_a \rangle$, we must have 
\begin{equation}
\label{eqn:F1}
  \w^u_0 = 0 \, ;
\end{equation}
and the forms $\w^0_0, \w^1_0, \ldots, \w^n_0$ are linearly independent.
Differentiating this formula and an application of Cartan's Lemma 
(cf. \cite{BCG3} or \cite[Lem. A.1.9]{IL03}) produces 
\begin{equation}
\label{eqn:F2}
  \w^u_a = q^u_{ab} \, \w^a_0 \, ,
\end{equation}
where, $q^u_{ab} = q^u_{ba}$ are $\C$-valued functions on $\F^1$.  The 
$q^u_{ab}$ are the coefficients of the {\it second fundamental form}
\begin{displaymath}
  II = q^u_{ab} \, \w^a_0 \w^b_0 \, \otimes \, \underline e_{u} 
  \in \Gamma ( \F^1 , \pi^*( S^2 \, T^{*}M \otimes NM ) )\, ,   
\end{displaymath}
a section of the pulled-back bundle $\pi^*(S^2 \, T^{*}M \otimes NM )$ over 
$\F^1$.  Here $TM$ and $NM$ denote the tangent and normal bundles over $M$.
Given $x = [e_0] \in M$, 
$T_xM \simeq \hat x^{*} \otimes (T_{e_0} \widehat M / \hat x)$; the normal 
space $N_xM = T_x\P V/T_xM$ is spanned by 
$\underline e_u := e^0 \otimes ( e_u \hbox{ mod } T_{e_0} \widehat M )$; and 
the $e^I \in V^*$ are dual to the $e_I \in V$.
A priori defined on $\F^1$, the second fundamental form descends to a 
well-defined section of $S^2(T^*M) \otimes NM$ over $M$.  

Similarly, differentiating (\ref{eqn:F2}) yields the {\it Fubini cubic form} 
\begin{displaymath}
  F_3 = r^u_{abc} \, \w^a_0 \w^b_0 \w^c_0 \, \otimes \underline e_u 
  \in \Gamma ( \F^1 , \pi^*( S^3 \, T^*M \otimes NM ) )\, ,
\end{displaymath}
a higher order differential
invariant.  The coefficients $r^u_{abc} : \F^1 \to \C$ of $F_3$ are fully 
symmetric in the lower indices, and are defined by 
\begin{equation}
\label{eqn:F3}
  r^u_{abc} \w^c_0 = - d q^u_{ab} - q^u_{ab} \, \w^0_0 - q^v_{ab} \, \w^u_v
                     + q^u_{ae} \, \w^e_b + q^u_{be} \, \w^e_a \, . 
\end{equation}
Continuing in this fashion we may construct a sequence of higher order 
differential invariants 
\begin{displaymath}
  F_k = r^u_{a_1 \cdots a_k} \, \w^{a_1}_0 \cdots \w^{a_k}_0 \, \otimes
        \underline e_u 
  \in \Gamma ( \F^1 , \pi^*( S^k \, T^*M \otimes NM ) ) \, , \ k \ge 3 \, . 
\end{displaymath}
%
%
\begin{prop}
\label{prop:F}
Set $p = 2q$ if $p$ is even, and $p = 2q + 1$ if 
$p$ is odd.  Then the coefficients of $F_{p+1}$, again fully
symmetric in the lower indices, are defined by
\begin{eqnarray*}
  r^u_{a_1 \ldots a_p b} \, \w^b_0 & = &
  - \, d r^u_{a_1 \ldots a_p} \, - \, (p-1) \, r^u_{a_1 \ldots a_p} \, \w^0_0
  \, - \, r^v_{a_1 \ldots a_p} \, \w^u_v \hfill \\
  & &
  + \ \fS_{a_1 \ldots a_p} 
  \left\{ (p-2) \, r^u_{a_1 \ldots a_{p-1}} \w^0_{a_p} \,  
       + \, r^u_{a_1 \ldots a_{p-1} b} \, \w^b_{a_p} \, 
       - \, q^u_{a_p b} \, r^v_{a_1 \ldots a_{p-1}} \, \w^b_v \right\} \\
  & & 
  - \sum_{k=q+1}^{p-2} \, \cS_{a_1 \ldots a_k , a_{k+1} \ldots a_p}
    \bigg\{ \hspace{5pt}
    \left( r^u_{a_1 \ldots a_k b} \, r^v_{a_{k+1} \ldots a_p} \, + \, 
           r^u_{a_{k+1} \ldots a_p b} \, r^v_{a_1 \ldots a_k} 
    \right) \, \w^b_v \\
  & & \hspace{65pt}
    + \ \Big(
             (k-1) \, r^u_{a_1 \ldots a_k} \, r^v_{a_{k+1} \ldots a_p} \
         + \ (p-k-1) \, r^u_{a_{k+1} \ldots a_p} \, r^v_{a_1 \ldots a_k} 
         \Big) \, \w^0_v
    \bigg\} \\
  & & - \ \tbinom{p+1}{mod \ 2} \, 
      \cS_{a_1 \ldots a_q , a_{q+1} \ldots a_p}
      \, \bigg\{ 
      \ r^u_{a_1 \ldots a_q b} \, r^v_{a_{q+1} \ldots a_p} \, \w^b_v \ 
      + \ (q-1)\, r^u_{a_1 \ldots a_q} \, r^v_{a_{q+1} \ldots a_p} \, 
                 \w^0_v \, \bigg\} \, . 
\end{eqnarray*}
\end{prop}
As in \cite{IL03} $\fS$ denotes cyclic summation over the indices.
The notation $\cS$ represents a symmetrizing operation defined as follows:
Given symmetric tensors $T_{a_1 \ldots a_k}$ and $U_{b_1 \ldots b_\ell}$, 
$\cS_{a_1 \ldots a_k,b_1 \ldots b_\ell}$ denotes summation over the 
$\tbinom{k+\ell}{k} = \tbinom{k+\ell}{\ell}$ 
elements of 
the symmetric group on $k+\ell$ elements symmetrizing the product.  For 
example, $\cS_{a_1 \ldots a_k , b} = \fS_{a_1 \ldots a_k b}$ and 
$\cS_{ab,cd} = \fS_{abc} \cdot \fS_{bd}$.  

Notice that the last line in the equation appears only when $p$ is even.  
Also, I've used the convention $II = F_2$ and $q^u_{ab} = r^u_{ab}$.
\vspace{0.1in}\\
\noindent
{\it Remark.}  This formula corrects errors in the expression for 
$r^u_{a_1 \ldots a_p b} \w^b_0$ given in \cite[p. 108]{IL03} and 
\cite[(2.20)]{L96}.

\begin{proof}
The formula is easily verified for $F_3$ and $F_4$.  (For $F_3$ it is necessary
to use the convention that $r^u_a = 0$.)  Obtain the coefficients of $F_4$ 
\begin{equation}
\label{eqn:F4}
  r^u_{abcd} \w^d_0 = -d r^u_{abc} - 2 r^u_{abc} \w^0_0 - r^v_{abc} \w^u_v
  + \fS_{abc} \left(
    r^u_{abe} \w^e_c + q^u_{ab} \w^0_c - q^u_{ae} q^v_{bc} \w^e_v 
  \right)
\end{equation}
by differentiating (\ref{eqn:F3}) and another application of Cartan's Lemma.  
The general statement is established by induction; once for $p$ even, and once 
for $p$ odd.
\end{proof}

\noindent
The formula for $F_k$ given by Proposition \ref{prop:F} provides a slight 
improvement of \cite[Prop. 2.41]{L96}:
%
%
\begin{cor}
\label{cor:F}
Let $x$ be a smooth point of a variety $X \subset \P V$.  Suppose that there 
exists a framing $\textbf{e}_x=(e_0,\ldots,e_N)$ over $x = [e_0]$ where the 
coefficients of $F_\ell, F_{\ell + 1}, \ldots, F_{2\ell - 1}$ all vanish.  
If $x$ is a $(2\ell-1)$-general point, then 
the coefficients of $F_k$ vanish on $\F^{2\ell-1}$, for all $k \ge \l$.
\end{cor}

\noindent
{\it Remark}.  Loosely speaking, after taking $k$ derivatives there will be 
both discrete and continuous invariants: the discrete invariants are locally
constant at a {\it $k$-general point}. 

\begin{proof}
This follows immediately from Proposition \ref{prop:F}.
\end{proof}

Unlike the second fundamental form, the 
$F_k \in \Gamma(\F^1 , \pi^*( S^k \, T^*X \, \otimes \, NX ))$, $k\ge3$, do 
not descend to well-defined sections over $X$.  For this reason we call them 
{\it relative differential invariants}.  However:
\begin{itemize}
  \item Let $|F_{x,\l}| = F_\l(N_x X) \subset S^\l( T^*_xX)$.  The zero locus 
    $\mathcal C_{k,x}$ of $|II_x| , |F_{3,x}| , \ldots , |F_{k,x}|$ in 
    $\P( T_xX)$ is well defined, and consists of the tangent directions to 
    lines making contact to order $k$ with $X$ at $x$.  See 
    \S\ref{sec:rigidity} below.
  \item By restricting $F_3 : N^*X \to S^3 \, T^*X$ to the kernel of 
    $F_2 = II : N^*X \to S^2 \, T^*X$ we obtain a tensor 
    $\mathbb{F}_3 = III \in S^3 \, T^*X \otimes N_3$ on $X$.  Here 
    $N_3 = T_x \P V / \{ T_xX \oplus II(S^2 \, T_xX) \}$, and $\mathbb F_3$
    is called the {\it third fundamental form}.  A series of higher order 
    fundamental forms $\mathbb F_k \in S^k \, T^*X \otimes N_k$ on $X$ is 
    defined inductively.  See \cite[\S3.5]{IL03} for details.    
\end{itemize}


\section{Griffiths--Harris rigidity}
\label{sec:rigidity}

\noindent
The $F_k$ play an important role in establishing the rigidity of a variety.
%
%
\begin{prop}[{\cite[Cor. 3.7.2]{IL03}}]
\label{prop:TS}
A complex projective variety $X \subset \P V$ is uniquely determined up to 
projective equivalence by the infinite sequence of relative differential 
invariants at a smooth point $x \in X$.
\end{prop}

\noindent
The proof of the proposition is straightforward.  Here is a sketch.  Near a 
smooth point $x$, $X$ may be expressed locally as a graph 
$x^u = f^u(x^1 , \ldots , x^n)$.  There exists a local section 
$s : X \to \F^1$ so that the pull-back of $F_k$ is 
$s^*(F_k) = (-1)^k \frac{\partial^k f^u}{\partial x^{a_1} \cdots 
                   \partial x^{a_k}} \, dx^{a_1} \cdots dx^{a_k} \, \otimes \, 
                   \frac{\partial}{\partial x^u}$ 
at $x$.  In particular, the coefficients of $F_k$ determine the Taylor 
series of the $f^u$, which in turn determine $X$.
\vspace{0.1in}

\begin{dfn}
  Two varieties $X^n, Y^n \subset \P^N$ agree to order $k$ at 
  $x\in X_\textit{smooth}$ and $y\in Y_\textit{smooth}$ if there exist frames 
  $\textbf{e}_x\in \F^1_X$ and $\textbf{e}_y\in \F^1_Y$ over 
  $x$ and $y$ such that the coefficients of $F_\l(\textbf{e}_x)$ and 
  $F_\l(\textbf{e}_y)$ are equal for all $\l \le k$.
\end{dfn}
\noindent
Notice this implies $\mathcal C_{\l,x} = \mathcal C_{\l,y}$, $\l \le k$.
\begin{dfn}
  When agreement to $k$-th order forces agreement to all orders, then 
  Proposition \ref{prop:TS} implies that $X$ and $Y$ are projectively 
  equivalent, and we say $X$ is rigid to order $k$.
\end{dfn}

\noindent
{\it Remark.}  Any variety $X$ meeting the conditions of Corollary \ref{cor:F} 
is rigid to order $k = 2\ell - 1$.
\vspace{0.1in}

Here is another perspective on rigidity.  All of the varieties $X$ to be 
discussed in this paper admit sub-bundles $\F_X \subset \F^1_X$ of the 
first-order adapted frames on which the coefficients of the $F_k$ are constant.
The sub-bundles $\F_X$ are maximal integral submanifolds of a Frobenius system 
$\mathcal I$ on $\mathrm{GL}(V)$.  The Frobenius system is generated by 
constant coefficient linear combinations of the entries of the Maurer-Cartan 
form on $\mathrm{GL}(V)$.  In particular, any other maximal integral 
submanifold of $\mathcal I$ is of the form $g\cdot \F_X$, $g\in\mathrm{GL}(V)$.

Suppose $Y^n$ agrees with $X^n$ (both algebraic varieties of $\P^N$) to order 
$k$ at a $k$-general point $y \subset Y_\textit{smooth}$.  Equivalently, we may
restrict to a sub-bundle $\F^k_U \subset \F^1_Y$ over an open neighborhood 
$U \subset Y_\textit{smooth}$ of $y$ on which the the coefficients of 
$F_{\l,Y}$ are equal to the constant coefficients of $F_{\l,X}$, $\l \le k$.  
The variety $X$ will be rigid to order $k$ if and only if we may further 
reduce (or {\it normalize}, see \S\ref{sec:normalizations}) $\F^k_U$ to a 
sub-bundle $\F_U$ which is an integral submanifold of $\mathcal I$.  When this 
is the case, we have $g \cdot \F_U \subset \F_X$ yielding a projective linear
transformation $g\cdot U \subset X$.  As we are working with algebraic
varieties, it now follows that $g \cdot Y = X$.

\subsection{Previous rigidity results}

\subsubsection{The Segre variety}

Let $W_1^*$ and $W_2$ be complex vector spaces of dimensions $d_1, d_2>1$, 
respectively.  The Segre variety 
$X = \textrm{Seg}(\P W_1^* \times \P W_2) = \{ [ w_1^* \otimes w_2] \} \subset 
\P(W_1 \otimes W_2)$ is the set 
$\P\{\hbox{rank one linear maps } W_1 \to W_2 \}$.
Given $d_2 > 2$, Landsberg has shown that $X$ is rigid to order 2 \cite{L99}.
For all $k \ge 3$, the $F_k$ may be normalized to zero; and 
$\mathcal C_x = \mathcal C_{2,x} = \P^{d_1-1} \sqcup \P^{d_2-1}$ is the 
disjoint union of two linear subspaces.

\subsubsection{The Veronese variety}
\label{sec:veronese}

Let $W$ be a complex vector space of dimension $n+1>2$, and consider the 
Veronese embedding $X = v_2(\P W) \subset \P(S^2\, W)$ of $\P W$.  In this 
case $|II_x| = S^2 \, T^{*}_xX$, so that $\mathcal C_{2,x} = \emptyset$, and 
there exists a sub-bundle of the first order adapted frames upon which the 
$F_k = 0$, $k\ge 3$.  Landsberg has shown that $v_2(\P W)$ is rigid to order 
three \cite{L99}. (The case $v_2(\P^2)\subset\P^5$ was established by Griffiths
and Harris \cite{GH79}.)  Note that, in the case that $n = 2m-1$, $X$ is the 
adjoint variety of of the simple Lie group $\textrm{Sp}(W)$.

\subsubsection{Compact Hermitian symmetric spaces}
\label{sec:cHss}

It is a classic result of Fubini that a quadric hypersurface (of rank $> 1$) 
in $\P^N$, $N > 2$, is rigid to order three; cf. \cite{F18} or 
\cite[Th. 3.9.1]{IL03}.  (In contrast, Monge showed that the conics in $\P^2$ 
are rigid to order five \cite[\S3.6]{L96}.)  
The quadric hypersurfaces are rank two compact Hermitian symmetric spaces.
Landsberg has shown that any other rank two CHSS, in its minimal homogeneous 
embedding, is rigid to order two \cite{L06}.

Shortly after Landsberg's preprint appeared Hwang and Yamaguchi proved the 
following theorem:  Let $Y^n \subset \C\P^{N}$ be a complex submanifold and 
$y\in Y$ a general point.  Suppose that $X$ is an irreducible rank $r$ compact
Hermitian symmetric space in its natural embedding, other than a quadric 
hypersurface.  If there exist linear maps $f:T_xX \to T_yY$ and 
$g_k: N_{k,x}X \to N_{x,y}Y$ such that the induced maps 
$S^k \, T^*_xX \otimes N_{k,x} X \to S^k \, T^*_yY \otimes N_{k,y}$ take 
$\mathbb F_{k,x} \to \mathbb F_{k,y}$, for $2 \le k \le r$, then 
$\overline Y$ is projectively equivalent to $X$ \cite{HY03}.  Unlike the other
results mentioned in this section, Hwang and Yamaguchi's proof is not based 
on a moving frame calculation.  Their result follows from an elegant 
application of representation theory, based on a theory developed by Se-ashi 
\cite{S88}, and a Lie algebra cohomology computation made easy by the work of
Kostant \cite{K61}.

The fundamental forms $\mathbb F_k$ are the only nonzero invariants of the 
CHSS.  In particular, the $F_k$, $k \ge 3$, may be normalized to zero, so that 
$\mathcal C_{x} = \mathcal C_{2,x}$.

The interested reader will find {\it intrinsic} rigidity results for the  
compact Hermitian symmetric spaces in the paper \cite{HM06} of Hwang and Mok.

\section{Normalizations}
\label{sec:normalizations}

Given a frame $\textbf{e}$ in the fibre over a smooth point $x=[e_0] \in X$,  
consider the fibre motion $\t{\textbf{e}} = g \cdot \textbf{e}$, $g\in GL(V)$. 
It will be helpful to compare the expressions for $F_k$ at $\textbf{e}$ and 
$\t{\textbf{e}}$.  Transformations by block diagonal matrices 
$\t e_0 = g^0_0 e_0$, $\t e_a = g^b_a e^b$, and $\t e_u = g^v_u e_v$ do not 
change the $F_k$. For example, if the coefficients of $II$ are $q^u_{ab}$ at 
$e$, then the coefficients of $II$ at $\t e$ are 
$\t q^u_{ab} = h^0_0 h^u_v g^c_a g^d_b q^v_{cd}$, and $II_e = II_{\t e}$.
So we consider transformations of the form
\begin{equation}
\label{eqn:tilde}
  \t e_0 = e_0 \, , \quad 
  \t e_a = e_a + g^0_a e_0 \, , \quad 
  \t e_u = e_u + g^0_u e_0 + g^a_u e_a \, .
\end{equation}
Let $\t \w$ denote the Maurer-Cartan forms at $\t e \in \F^1$, so that 
$d \t e = \t \w \t e$.  Making use of $\t e = g \, e$, we derive 
$\t \w = g^{-1} \w  g + g^{-1} \, dg$.  Explicitly,
\begin{equation}
\label{eqn:tildew}
\renewcommand{\arraystretch}{1.5}
\begin{array}{rp{20pt}l}
  \t \w^a_0 = \w^a_0  & & 
  \t \w^u_b = \w^u_b  \\
  \t \w^0_0 = \w^0_0 - g^0_c \, \w^c_0 & & 
  \t \w^a_b = \w^a_b + g^0_b \, \w^a_0 - g^a_t \, q^t_{bc} \, \w^c_0 \\ 
  \t \w^u_v = \w^u_v + g^e_v \, q^u_{ce} \, \w^c_0 & &
  \t \w^a_v =\w^a_v + ( \d^a_c \, g^0_v  - g^e_v \, g^a_t \, q^t_{ce} ) \, 
                      \w^c_0 + g^c_v \, \w^a_c - g^a_t \, \w^t_v 
                    + d g^a_v \\ 
  \multicolumn{3}{l}{ 
  \t \w^0_b = \w^0_b 
            + ( q^t_{bc} \, ( g^0_e \, g^e_t - g^0_t ) - g^0_b \, g^0_c ) \,
              \w^c_0 + g^0_b \, \w^0_0 - g^0_c \, \w^c_b 
            + d g^0_b } \\
  \multicolumn{3}{l}{
  \t \w^0_v
            = \w^0_v 
            + \left( 
              g^e_v \, ( g^0_d \, g^d_t - g^0_t ) \, q^t_{ce} - g^0_v \, g^0_c 
              \right) \, \w^c_0 
            + g^0_v \, \w^0_0 + g^c_v \, \w^0_c - g^0_c \, ( g^e_v \, \w^c_e 
            + \w^c_v ) } \\
  \multicolumn{3}{l}{\hspace{38pt} + \left( g^0_e \, g^e_t - g^0_t \right) \, 
              \w^t_v + d g^0_v - g^0_e \, d g^e_v }
\end{array}
\end{equation}
Using the two expressions on the first line we immediately see that 
\begin{equation}
\label{eqn:tilde2}
  \t q^u_{ab} = q^u_{ab} \, .
\end{equation}

The coefficients of $F_3$ at $\t e \in \F^1$ are given by (\ref{eqn:F3}) 
\begin{eqnarray*}
  \t r^u_{abc} \t\w^c_0  
  =  - d \t q^u_{ab} 
                            - \t q^u_{ab} \t\w^0_0 - \t q^v_{ab} \t \w^u_v
  + \t q^u_{ae} \t \w^e_b + \t q^u_{be} \t \w^e_a \, .
\end{eqnarray*}
Replace the $\t q$ and $\t\w$ with their $q$ and $\w$ expressions.  After 
simplifying, both sides of the equation are seen to be linear combinations of 
the $\w^c_0$.  Equating coefficients produces
\begin{equation}
\label{eqn:tilde3}
 \t r^u_{abc} = r^u_{abc} 
              + \mathfrak{S}_{abc} 
                ( \d^u_v g^0_a - q^u_{ae} g^e_v ) q^v_{bc} \, .
\end{equation}
Consequently, we see that although the subspace $|F_3| \subset S^3 \, T_x^*X$
is not well-defined over $x$ (the subspace moves as we vary the frame $e$ over
$x$), it is well-defined over $x$ modulo $|II| \circ T^*X$.  In particular, 
the zero locus $\mathcal C_{3,x} \subset \P(T_xX)$ of $\{ |II| , |F_3| \}$ is 
well-defined.

Similarly (\ref{eqn:F4},\ref{eqn:tildew}) yield,
\begin{displaymath}
\renewcommand{\arraystretch}{1.5}
\begin{array}{rcl}
  \t r^u_{abcd} & = & 
  r^u_{abcd} \ 
  + \ \mathfrak S_{abcd} \, 
    \Big\{ 
        2 \, r^u_{abc} g^0_d \,
      - \, r^v_{abc} \, q^u_{de} \, g^e_v \Big\} \\
  & & 
  + \ \mathcal{S}_{a,b,cd} \Big\{ \, 
             ( \d^u_v \, g^0_a \, - \, q^u_{ae} \, g^e_v ) \, 
             ( \d^v_w \, g^0_b \, - \, q^v_{be} \, g^e_w ) \, 
             q^w_{cd}
                              \Big\} \\
  & & 
  + \ \mathfrak{S}_{abc} \, \Big\{ \ 
      q^u_{e_1e_2} \, g^{e_1}_v \, g^{e_2}_w \, q^v_{ab} \, q^w_{cd} \ 
                             \Big\}
  \ - \
       \mathcal{S}_{ab,cd} \, \Big\{
       ( q^u_{ab} \, g^0_v \, + \, r^u_{abe} \, g^e_v ) \,  q^v_{cd}  
                               \Big\} \, , \\
\end{array}
\end{displaymath}
with $\mathcal{S}_{a,b,cd} = \mathfrak{S}_{abcd} \cdot \mathfrak{S}_{bcd}$.
This expression corrects typos in the formula for $\t r^u_{abcd}$ in 
\cite[p. 108]{IL03}.  As in the case of $F_3$, the subspace 
$|F_4| \subset S^4 \, T_x^* X$ is well-defined modulo the ideal generated by 
$|II|, |F_3|$.

In this paper, we will be interested in the {\it special case} that 
$g^0_b = 0 = g^a_v$; from (\ref{eqn:tilde2},\ ref{eqn:tilde3}) we see that 
these are the fibre motions that preserve the coefficients of $II$ and $F_3$.  
Then the fibre variation for the coefficients of $F_4$ is given by 
\begin{equation}
\label{eqn:tilde4}
  \t r^u_{abcd} \ = \  
  r^u_{abcd} \ - \ 
  \mathcal{S}_{ab,cd} \, \left( q^u_{ab} \, g^0_v \, q^v_{cd}  \right) \, .
\end{equation}

\section{The trace-free, rank 1 matrices}
\label{sec:SL}

Let $W$ be a complex vector space of dimension $m+1$, and let $X$ denote the 
variety of trace-free, rank one linear maps $W \to W$.  Observe that the rank 
one transformations may be identified with the Segre variety
$\hbox{Seg}(\P W \times \P W^*) = \{ [v\otimes w^*] \ | \ v \in W \, , \ 
w^* \in W^* \}$.  The matrix $v \otimes w^*$ is trace-free if and only if 
$w^*(v) = 0$.  That is, $X$ is a hyperplane section of the Segre variety.

This variety is the unique closed orbit of the adjoint action on $\gsl$.
Let $\{ v_j \}_{j=0}^m$ be a basis of $W$ and $\{ v^*_j \}$ the dual basis.
Let $f = (f_0 , \ldots , f_m) \in \GSL$, and write 
$f^{-1} = (f_0^* , \ldots , f_m^*)^t$.  Then the orbit of 
$v_0 \otimes v^*_m \in \gsl$ under the adjoint action is 
\begin{center}
  $\{ f_0 \otimes f^*_m \ | \ 
  f \in \GSL \} 
  \subset \gsl$; 
\end{center}
and this is precisely the set of trace-free, rank one matrices.

Fix the index ranges
\begin{displaymath}
\renewcommand{\arraystretch}{1.4}
\begin{array}{rcccl}
  0 & \le &  j , k & \le & m \, , \\
  1 & \le & \a , \b & \le & m-1 \ , \\
  m & \le & \oa , \ob & \le & 2m-2 \, ,
  \quad \oa = \a + (m-1) \, , \\
  1 & \le & a , b & \le & 2m-1 = n \, .
\end{array}
\end{displaymath}
The vectors 
\begin{equation}
\label{eqn:SLframe}
\renewcommand{\arraystretch}{1.5}
\begin{array}{lll}
  e_0      = f_0 \otimes f_m^*  \, ;                                  & & \\
  e_\a     = f_\a \otimes f_m^* \, ,                                  &
  e_\oa    = f_0 \otimes f_\a^* \, ,                                  &
  e_{2m-1} = \half \, ( f_m \otimes f_m^* - f_0 \otimes f_0^* ) \, ;  \\  
  e_{\a0}  = f_\a \otimes f_0^* \, ,                                  &
  e_{m\a}  = f_m \otimes f_\a^* \, ,                                  &
  e_{m0}   = f_m \otimes f_0^*  \, ,                                  \\
  \multicolumn{3}{l}{
  e_{\a\b} = f_\a \otimes f_\b^* - \half \d_{\a\b} 
                 ( f_0 \otimes f_0^* + f_m \otimes f_m^* ) }
\end{array}
\end{equation}
span $\gsl$, yielding a map 
$\varphi : \textrm{SL}(W) \to \textrm{GL}(\mathfrak{sl}(W))$.

Let $\F \subset \textrm{GL}(\mathfrak{sl}_{m+1})$ denote the set 
of all such framings.  Observe that $\F$ is a sub-bundle of the first-order 
adapted framings over $X$.  To see this, let $\eta$ be the Maurer-Cartan form 
on $\GSL$, so that
\begin{displaymath}
  d f_j = \eta^k_j \, f_k \, , \quad \hbox{ and } \quad 
  d f_j^* = -\eta^j_k \, f_k^* \, .
\end{displaymath}
These 1-forms satisfy $\eta^j_j = 0$, and are otherwise linearly independent.
Notice that $e_0 \in \widehat X$, and 
\begin{displaymath}
  d e_0 = \eta_0^j f_j \otimes f_m^* - f_0 \otimes \eta_j^m f_j^*
        = \left( \eta^0_0 - \eta^m_m \right) e_0
          + \eta_0^\a \, e_\a - \eta^m_\a \, e_\oa 
          + 2 \, \eta_0^m \,  e_{2m-1} \, .
\end{displaymath}
This implies that the $e_\a$, $e_\oa$, $e_{2m-1}$ span $T_{e_0} \widehat X$, 
and we have a first-order adapted framing of $\gsl$ over $X$.


\subsection{The Maurer-Cartan form on \boldmath $\F$ \unboldmath}

Now let $\w$ denote the pull-back of the Maurer-Cartan form on 
$\hbox{GL}(\mathfrak{sl}_{m+1})$ to $\F$, so that $d e = \w e$.  Computations 
analogous to that of $d e_0$ above assure us that the 1-forms $\w^0_0$, 
$\w^a_0$, $\w^0_a$ and $\w^\a_\b$ are linearly independent, and yield the 
following relations
\begin{equation}
\label{eqn:SLMC}
\renewcommand{\arraystretch}{1.5}
\begin{array}{l}
  \w^{2m-1}_\a = -\w^\oa_0 \, , \quad
  \w^{\b0}_\a = - \half \, \d^\b_\a \, \w^{2m-1}_0 \, , \quad
  \w^{\b\c}_\a = \d^\b_\a \, \w^\oc_0   
  \\
  \w^{2m-1}_\oa = \w^\a_0 \, , \quad
  \w^{m\b}_\oa = \half \, \d^\b_\a \, \w^{2m-1}_0 \, , \quad
  \w^{\b\c}_{\oa} = \d^\c_\a \, \w^\b_0 \, , \quad
  \w^\ob_\oa = -\w^\a_\b + \d_{\a\b} \, \w^0_0 
  \\
  \w^\a_{2m-1} = -\half \, \w^0_\oa \, , \quad
  \w^\oa_{2m-1} = \half \, \w^0_\a \, , \quad
  \w^{\a 0}_{2m-1} = -\half \, \w^\a_0 \, , \quad
  \w^{m\a}_{2m-1} = \half \, \w^\oa_0 \, , \quad
  \w^{m0}_{2m-1} = -\half \, \w^{2m-1}_0 
  \\
  \w^\b_{\a0} = -\d^\b_\a \, \w^0_{2m-1} \, , \quad 
  \w^{2m-1}_{\a0} = -\w^0_\a \, , \quad
  \w^{\b0}_{\a0} = \w^\b_\a - \d^\b_\a \, \w^0_0 \, , \quad
  \w^{\b\c}_{\a0} = -\d^\b_\a \, \w^0_\c \, , \quad 
  \w^{m0}_{\a0} = -\w^\oa_0 
  \\
  \w^\ob_{m\a} = \d^\b_\a \, \w^0_{2m-1} \, , \quad 
  \w^{2m-1}_{m\a} = \w^0_\oa \, , \quad
  \w^{m\b}_{m\a} = -\w^\a_\b \, , \quad
  \w^{\b\c}_{m\a} = -\d^\c_\a \, \w^0_\ob \, , \quad 
  \w^{m0}_{m\a} = -\w^\a_0 
  \\
  \w^\c_{\a\b} = \left( \d^\c_\a \, \d^\e_\b 
    + \half \, \d_{\a\b} \, \d^{\c\e} \right) \, \w^0_\oe \, , \quad 
  \w^\oc_{\a\b} = \left( \d^\e_\a \, \d^\c_\b 
    + \half \, \d_{\a\b} \, \d^{\c\e} \right) \, \w^0_\e
  \\
  \w^{\c0}_{\a\b} = - \left( \d^\c_\a \, \d_{\b\e} 
    + \half \, \d_{\a\b} \, \d^\c_\e \right) \, \w^\e_0 \, , \quad
  \w^{m\c}_{\a\b} = - \left( \d_{\a\e} \, \d^\c_\b 
    + \half \, \d_{\a\b} \, \d^\c_\e \right) \, \w^\oe_0 \, , \quad 
  \w^{\c\e}_{\a\b} = \d^\e_\b \, \w^\c_\a - \d^\c_\a \, \w^\b_\e
  \\
  \w^{2m-1}_{m0} = -2 \, \w^0_{2m-1} \, , \quad 
  \w^{\a0}_{m0} = -\w^0_\oa \, , \quad
  \w^{m\a}_{m0} = -\w^0_\a \, , \quad 
  \w^{m0}_{m0} = -\w^0_0
  \end{array}
\end{equation}
The remaining 1-forms (those not appearing above) vanish on the pull-back.

These are precisely the equations of the $\varphi$-pullback of the 
Maurer-Cartan form on $\textrm{GL}(\mathfrak{sl}(W))$.  In fact, if 
$\eta = \eta^{j}_{k} E_{j}^{k}$ denotes the Maurer-Cartan form on 
$\textrm{SL}(W)$, then
\begin{center}
\renewcommand{\arraystretch}{1.8}
\begin{tabular}{p{1.25in}p{1.25in}p{1.25in}p{1.25in}}
  \multicolumn{2}{l}{$\displaystyle
  \eta^{0}_{0} = \half \, \left(
      \textstyle{\frac{1}{n+1}} \sum \left( \w^{\oa}_{\oa} - \w^{\a}_{\a} \right) + \w^{0}_{0}
    \right)$ }  & 
  $\eta^{0}_{\b} = \w^0_\b$  & 
  $\eta^{0}_{n} = \w^0_{2m-1}$ 
  \\ 
  $\eta^{\a}_{0} = \w^\a_0$  &
  \multicolumn{2}{l}{$\displaystyle
    \eta^{\a}_{\b} = \half \, \left(
      \w^{\a}_{\b} - \w^{\ob}_{\oa} + \textstyle{\frac{\d^{\a\b}}{n+1}} \sum \left( \w^{\oa}_{\oa} - \w^{\a}_{\a} \right)
    \right)$ } &
  $\eta^{\a}_{n} = -\w^0_\oa$
  \\ 
  $\eta^{n}_{0} = \half \w^{2m-1}_0$ & 
  $\eta^{n}_{\b} = -\w^\ob_0$ &
  \multicolumn{2}{l}{$\displaystyle 
  \eta^{n}_{n} =  \half \, \left(
      \textstyle{\frac{1}{n+1}} \sum \left( \w^{\oa}_{\oa} - \w^{\a}_{\a} \right) - \w^{0}_{0}
    \right) \, .$}
\end{tabular}
\end{center}


\subsection{The differential invariants \boldmath $F_k$\unboldmath}
\label{sec:SLFk}

Recollect (\ref{eqn:F2}) that the coefficients $q^u_{ab}$ of the second 
fundamental form are defined by $\w^u_a = q^u_{ab} \w^b_0$.  Inspecting 
(\ref{eqn:SLMC}), we see that the non-zero coefficients are 
\begin{equation}
\label{eqn:SLF2}
  q^{m0}_{2m-1,2m-1} = -\half \ , \quad
  q^{\a0}_{\b,2m-1} = -\half \, \d^\a_\b \ , \quad
  q^{m\a}_{\ob,2m-1} = \half \, \d^\a_\b \ , \quad 
  q^{\a\b}_{\c\oe} = \d^\a_\c \, \d^\b_\e \, .
\end{equation}
In particular, 
\begin{displaymath}
  |II| = \langle 
         \w^{2m-1}_0 \, \w^{2m-1}_0 \, , \ 
         \w^\a_0 \, \w^{2m-1}_0 \, , \ 
         \w^\oa_0 \, \w^{2m-1}_0 \, , \ 
         \w^\a_0 \, \w^\oa_0 \rangle \subset S^2 \, T_xX \, .
\end{displaymath}
The first quadric in $|II|$ above implies that the cone over 
$\mathcal C_{2,x} = \textrm{Baseloc}|II|$ lies in the {\it contact hyperplane}
$T_1 := \langle \w^{2m-1}_0 \rangle^\perp \subset T_xX$.  The fourth quadric
tells us that $\mathcal C_{2,x}$, the set of lines 
osculating to order two at $x \in X$, is the disjoint 
union of two linear spaces
\begin{displaymath}
  \mathcal C_{2,x} = \P^{m-1} \sqcup \P^{m-1} \subset 
  \langle \w^{2m-1}_0 \rangle^\perp \subset T_xX \, . 
\end{displaymath}
It is straightforward to confirm that the two disjoint $\P^{m-1}$'s making up 
$\mathcal C_{2,x}$ correspond to integrable distributions, 
$D_1 = \{ 0 = \w^\a_0 \}$ and $D_2 = \{ 0 = \w^\oa_0 \}$, in $T_1$.

Computations with (\ref{eqn:F3}, \ref{eqn:F4}) show that the non-zero 
coefficients of $F_3$ are 
\begin{equation}
\label{eqn:SLF3}
  r^{\a 0}_{\b\c\oe} = \half ( \d^\a_\b \d_{\c\e} + \d^\a_\c \d_{\b\e} )
  \quad \hbox{ and } \quad
  r^{m\a}_{\b \oc \oe} = \half ( \d^\a_\c \d_{\b\e} + \d^\a_\e \d_{\b\c} ) \, ;
\end{equation}
and the nonzero coefficients of $F_4$ are 
\begin{equation}
\label{eqn:SLF4}
  r^{m0}_{\a\b\oc\oe} = \half ( \d_{\a\c} \d_{\b\e} + \d_{\a\e} \d_{\b\c} )
  = r^{\e 0}_{\a\b\oc} = r^{m\e}_{\oa\ob\c} \, .
\end{equation}
Therefore, 
$|F_{3,x}| = \langle \w^\a_0 \, (\sum \w^\e_0 \, \w^\oe_0 ) \, , \, 
                     \w^\oa_0 \, (\sum \w^\e_0 \, \w^\oe_0 ) \rangle$ and 
$|F_{4,x}| = \langle (\sum \w^\e_0 \, \w^\oe_0 )^2 \rangle$.  Whence 
$\mathcal C_{4,x} = \mathcal C_{3,x} = \mathcal C_{2,x}$.  Notice that cubics 
of $|F_{3,x}|$ are the derivatives of the $F_{4,x}$ quartic.

Finally, $F_k = 0$ for all $k \ge 5$.  (Compute $F_k = 0$ directly for 
$5 \le k \le 9$, and then apply Corollary \ref{cor:F}.)  Hence 
$\mathcal C_{x} = \mathcal C_{2,x}$.  That is, any line osculating to order 
two at $x\in X$ is necessarily contained in $X$.  This is consistent with the 
fact that $X$ is generated by degree two polynomials: the rank one matrices 
are given by the vanishing of their 2-by-2 minors.
\medskip

\noindent
The trace-free matrices are rigid to order three:

\begin{thm}
\label{thm:SL}
Let $Y^{2m-1} \subset \P(\mathfrak{sl}_{m+1}) = \P^{m^2+2m-1}$ be an algebraic
variety, with $m>1$.  Suppose that the exists a framing $\textbf{e}_y$ over a 
$3$-general point $y \subset Y$ at which the nonzero coefficients of 
$II_Y$ and $F_{3,Y}$ are given by equations (\ref{eqn:SLF2}) and 
(\ref{eqn:SLF3}), respectively.  Then $Y$ is projectively equivalent to the 
variety of trace-free, rank one $(m+1) \times (m+1)$ matrices.
\end{thm}

\noindent
{\it Remark.}  This result is better than had been expected.  Notice that the 
coefficients of $F_4$ can not be normalized to zero.  This non-vanishing led 
Landsberg and Manivel to conjecture that the the adjoint variety of 
$\mathrm{SL}_{m+1}\C$ was rigid to order four, but not to order three
\cite{LM02}.  Indeed, this is the first example of a variety that is rigid to 
order $k$ for which the higher order $F_\l$, $\l > k$, {\it can not} all be 
normalized to zero.
\vspace{0.1in}\\

\noindent
{\it Remark.}  It suffices to assume that the coefficients of $II$ and $F_3$
may be but in the form (\ref{eqn:SLF2},\ref{eqn:SLF3}) at a frame $e\in\F^1_Y$
over a general point $y\in Y$.  (More precisely, $y$ is a 3-general point: the 
discrete invariants associated to second and third order data should be 
constant in a neighborhood of $y$.)
\vspace{0.1in}\\

\noindent
{\it Remark.}  Third order rigidity does not hold when $m=1$.  In this case we
have 
\begin{displaymath}
  f = \left( \begin{array}{rr}
        f^0_0 & f^0_1 \\
        f^1_0 & f^1_1 
      \end{array} \right) 
  \quad \hbox{ and } \quad 
  f^{-1} = \left( \begin{array}{rr}
        f^1_1 & -f^0_1 \\
        -f^1_0 & f^0_0 
      \end{array} \right) \, ;
\end{displaymath}
and the orbit of $v_0 \otimes v_1^*$ under the adjoint action is 
\begin{displaymath}
  \left\{ f_0 \otimes f_1^* = \left( \begin{array}{rr}
        -f^0_0 \, f^1_0 & f^0_0 \, f^0_0 \\
        -f^1_0 \, f^1_0 & f^1_0 \, f^0_0 
      \end{array} \right) \ : \ 
      f_0 = \left( \begin{array}{c} f_0^0 \\ f_0^1 \end{array} \right) 
      \in \C^2 \backslash \{0\} 
   \right\} \, .
\end{displaymath}
Notice that $f_0 \otimes f_1^*$ may be identified with the symmetric product 
$f_0 \circ f_0$.  Therefore the adjoint variety of $\textrm{SL}_2\C$ is the 
Veronese embedding $v_2(\P^1) \subset \P\, S^2 C^2 = \P^2$.  The plane conics 
are rigid to order five (cf. \S\ref{sec:cHss}).

\section{The proof of Theorem \ref{thm:SL}}
\label{sec:proof}

The goal of this section is to show that the adjoint variety
$X \subset \P(\mathfrak{sl}_{m+1})$ is rigid to order three.  That 
is, if $Y^{2m-1} \subset \P^{m^2+2m-1}$ 
admits a sub-bundle $\F^3_U$ over an open neighborhood $U\subset Y$ of $y$ of 
the first-order adapted frame bundle on which $II_Y = II_X$ and 
$F_{3,Y} = F_{3,X}$, then $Y$ is projectively equivalent to $X$.  Our strategy 
is to reduce $\F^3_U$ to a sub-bundle $\F^4_U$ on which
(i) the non-zero coefficients of $F_{4,Y}$ are given by (\ref{eqn:SLF4}), and 
(ii) the coefficients of $F_{k,Y}$ vanish for $5 \le k \le 9$.  Then Corollary
\ref{cor:F} and Proposition \ref{prop:TS} yield the desired rigidity.

The proof is structured as follows.  After deriving the consequences of 
third-order agreement, 
\begin{itemize}
  \item
    The fourth-order coefficients $r^{m0}_{abce}$ are computed in 
    \S\ref{sec:SLF4n0}.  Then all coefficients (except $r^{m0}_{\a\b\oc\oe}$, 
    which is given by (\ref{eqn:SLF4})) are normalized to zero.  We restrict 
    to the sub-bundle $\F^{4}_{Y} \subset \F^{3}_{Y}$ on which these 
    normalizations hold.
  \item 
    We show in \S\ref{sec:SLF4rz} and \S\ref{sec:SLF4rmdr} that the 
    coefficients $r^{\r\z}_{abce}$ and $r^{\r0}_{abce}, r^{m\z}_{abce}$ vanish 
    on $\F^{4}_{Y}$, respectively.
  \item
    In \S\ref{sec:SLF5+} we see that the coefficients of $F_{k,Y}$, $k \ge 5$ are zero on $\F^{4}_{Y}$.
  \item
    The calculations in \S\ref{sec:SLF4rz}--\ref{sec:SLF5+} assume that 
    $m > 2$.   I address the case $m=2$ in \S\ref{sec:SLsmalln}.
\end{itemize}
Let $\F^3_U$ be the sub-bundle of the first-order adapted frames on which 
the coefficients of $II_Y$ and $F_{3,Y}$ are given by (\ref{eqn:SLF2}, 
\ref{eqn:SLF3}).  As before, let $\w$ denote the pull-back of the 
Maurer-Cartan form on $\hbox{GL}(\mathfrak{sl}_{m+1})$ to $\F^3_U$.  The 
condition $II_Y = II_X$ implies
\begin{equation}
\label{eqn:SLorder2}
\renewcommand{\arraystretch}{1.5}
\begin{array}{l@{\hspace{20pt}}l@{\hspace{20pt}}l@{\hspace{20pt}}l}
  \w^{\b0}_\a = -\half \, \d^\b_\a \, \w^{2m-1}_0 & 
  \w^{m\b}_\a = 0 & 
  \w^{\b\c}_\a = \d^\b_\a \, \w^\oc_0 & 
  \w^{m0}_\a = 0 \\
  \w^{\b0}_\oa = 0 & 
  \w^{m\b}_\oa = \half \, \d^\b_\a \, \w^{2m-1}_0 & 
  \w^{\b\c}_\oa = \d^\c_\a \, \w^\b_0 & 
  \w^{m0}_\oa = 0 \\
  \w^{\b0}_{2m-1} = -\half \, \w^\b_0 &
  \w^{m\b}_{2m-1} = \half \, \w^\ob_0 &
  \w^{\b\c}_{2m-1} = 0 & 
  \w^{m0}_{2m-1} = -\half \, \w^{2m-1}_0 \, .
\end{array}
\end{equation}
The condition $F_{3,Y} = F_{3,X}$ yields
\begin{equation}
\label{eqn:SLorder3}
\renewcommand{\arraystretch}{1.5}
\begin{array}{l}
  0 \ = \ \w^\a_\ob \ = \ \w^\oa_\b \ = \
  \w^{\a0}_{m\b} \ = \ \w^{m\a}_{\b0} \ = \
  \w^{\a\b}_{m0} \ = \ \w^{m0}_{\a\b}
  \\
  \w^{2m-1}_\a \ = \ \w^{m0}_{\a0} \ = \ -\w^\oa_0 \, , \quad
  \w^{2m-1}_\oa \ = \ - \w^{m0}_{m\a} \ = \ \w^\a_0 \, , \quad 
  \w^0_0 + \w^{m0}_{m0} = 2 \, \w^{2m-1}_{2m-1}
  \\
  \w^{\a0}_{\b0} - \w^\a_\b = \d^\a_\b \, 
    \left( \w^{2m-1}_{2m-1} - \w^0_0 \right) \, , \quad
  \w^{\a0}_{\b\c} = - \left( \d^\a_\b \, \d_{\c\e} 
    + \half \, \d^\a_\e \, \d_{\b\c} \right) \, \w^\e_0 \, , \quad
  \w^{\a0}_{m0} = 2 \, \w^\a_{2m-1} 
  \\
  \w^{m\a}_{m\b} - \w^\oa_\ob = \d^\a_\b \, 
    \left( \w^{2m-1}_{2m-1} - \w^0_0 \right) \, , \quad
  \w^{m\a}_{\b\c} = - \left( \d^\a_\c \, \d_{\b\e} 
    + \half \, \d^\a_\e \, \d_{\b\c} \right) \, \w^\oe_0 \, , \quad
  \w^{m\a}_{m0} = -2 \, \w^\oa_{2m-1}
  \\
  \w^{\a\b}_{\c0} = -2 \, \d^\a_\c \, \w^\ob_{2m-1} \, , \quad 
  \w^{\a\b}_{m\c} = 2 \, \d^\b_\c \, \w^\a_{2m-1} \, , \quad
  \w^{\a\b}_{\c\e} + \d^\a_\c \, \d^\b_\e \, \w^0_0 = 
    \d^\a_\c \, \w^\ob_\oe + \d^\b_\e \, \w^\a_\c \, . 
\end{array}
\end{equation}
I will use these relations without mention when computing the coefficients of 
$F_{4,Y}$ below.


\subsection{$F_{4,Y}$ -- the conormal direction $u = m0$}  
\label{sec:SLF4n0}

Direct computations with (\ref{eqn:F4}) yield
\begin{eqnarray*}
  0 & = & r^{m0}_{\a\b\c e} \, \w^e_0 \ = \ 
  r^{m0}_{\oa\ob\oc e} \, \w^e_0 \ = \   
  r^{m0}_{\a\b(2m-1)e} \, \w^e_0 \ = \ 
  r^{m0}_{\oa\ob(2m-1)e} \, \w^e_0
  \\
  r^{m0}_{\a\b\oc e} \, \w^e_0 & = & 
    \half \, \left( \d_{\a\c} \, \d_{\b\e} + \d_{\a\e} \, \d_{\b\c} \right) 
    \w^\oe_0 
  \\
  r^{m0}_{\a\ob\oc e} \, \w^e_0 & = & 
    \half \, \left( \d_{\a\b} \, \d_{\c\e} + \d_{\a\c} \, \d_{\b\e} \right) 
    \w^\e_0 
  \\
  r^{m0}_{\a\ob(2m-1)e} \, \w^e_0 & = & \half \, \w^{2m-1}_{\a\b}
  \\
  r^{m0}_{\a(2m-1)^2e} \, \w^e_0 & = & - \half \, 
    \left( \w^{2m-1}_{\a0} + \w^0_\a \right)
  \\
  r^{m0}_{\oa(2m-1)^2e} \, \w^e_0 & = & \half \, 
    \left( \w^{2m-1}_{m\a} - \w^0_\oa \right) 
  \\ 
  r^{m0}_{(2m-1)^3 e} \, \w^e_0 & = & -\threehalves \, 
    \left( \w^0_{2m-1} + \half \w^{2m-1}_{m0} \right) \, .
\end{eqnarray*}
Equation (\ref{eqn:tilde4}) permits us to normalize the coefficients 
of the last four equations above to zero through transformations of the form 
(\ref{eqn:tilde}) with $g_a = 0 = g^a_u$ and 
\begin{displaymath}
  g^0_{\a0} = \textstyle{\frac43} r^{m0}_{(2m-1)^3\a} \ , \quad 
  g^0_{m\a} = -\textstyle{\frac43} r^{m0}_{(2m-1)^3\oa} \ , \quad 
  g^0_{\a\b} = -2 r^{m0}_{\a\ob(2m-1)^2} \ , \quad
  g^0_{m0} = \textstyle{\frac23} r^{m0}_{(2m-1)^4} \, .
\end{displaymath}
(Note that equations (\ref{eqn:tilde2}, \ref{eqn:tilde3}) assure us that the 
coefficients of $II_Y$ and $F_{3,Y}$ are preserved.)
Under these normalizations
\begin{equation}
\label{eqn:SL4n0}
  \w^{2m-1}_{\a0} = - \w^0_\a \ , \quad 
  \w^{2m-1}_{m\a} = \w^0_\oa \ , \quad 
  \w^{2m-1}_{\a\b} = 0 \ , \quad
  \w^{2m-1}_{m0} = -2 \w^0_{2m-1} \, ;
\end{equation}
and the only non-zero 
coefficients of $F_{4,Y}$ in the conormal direction $u = m0$ are given by (\ref{eqn:SLF4})
\begin{displaymath}
  r^{m0}_{\a\b\oc\oe} = \half \, 
    \left( \d_{\a\c} \, \d_{\b\e} + \d_{\a\e} \, \d_{\b\c} \right) \, .
\end{displaymath}
Restrict, from this point on, to the sub-bundle $\F^4_U \subset \F^3_U$ on which these normalizations hold.

\subsection{$F_{4,Y}$ -- the conormal direction $u = \r\z$}
\label{sec:SLF4rz}

After simplifying with (\ref{eqn:SLorder3}, \ref{eqn:SL4n0}), the coefficients $r^{\r\z}_{abce}$ are given by (\ref{eqn:F4}) as
\begin{eqnarray}
  0 & = & r^{\r\z}_{(2m-1)^3e} \, \w^e_0 \ = \ 
    r^{\r\z}_{\a\b\c e} \, \w^e_0  \ = \ 
    r^{\r\z}_{\oa\ob\oc e} \, \w^e_0  
  \label{eqn:sl1} \\
  r^{\r\z}_{\a(2m-1)^2e} \, \w^e_0 & = & \half \, \d^\r_\a \, \w^\oz_{m0} 
  \label{eqn:sl2} \\
  r^{\r\z}_{\oa(2m-1)^2e} \, \w^e_0 & = & \half \, \d^\z_\a \, \w^\r_{m0} 
  \label{eqn:sl3}  \\
  r^{\r\z}_{\a\b(2m-1)e} \, \w^e_0 & = & \half \, 
    \left( \d^\r_\a \, \w^\oz_{\b0} + \d^\r_\b \, \w^\oz_{\a0} \right) 
  \label{eqn:sl4} \\
  r^{\r\z}_{\a\ob(2m-1)e} \, \w^e_{0} & = & \d^\r_\a \, \d^\z_\b \, \w^0_{2m-1}
    + \half \, \d^\z_\b \, \w^\r_{\a0} - \half \, \d^\r_\a \, \w^\oz_{m\b}  
  \label{eqn:sl5} \\
  r^{\r\z}_{\oa\ob(2m-1)e} \, \w^e_0 & = & - \half \, 
    \left( \d^\z_\a \, \w^\r_{m\b} + \d^\z_\b \, \w^\r_{m\a} \right) 
  \label{eqn:sl6}  \\
  r^{\r\z}_{\a\b\oc e} \, \w^e_0 & = & 
    \d^\r_\a \, \left( \d^\z_\c \, \w^0_\b + \d_{\b\c} \, \w^\oz_{2m-1} 
      - \w^\oz_{\b\c} \right) +
    \d^\r_\b \, \left( \d^\z_\c \, \w^0_\a + \d_{\a\c} \, \w^\oz_{2m-1} 
      - \w^\oz_{\a\c} \right)  
  \label{eqn:sl7} \\
  r^{\r\z}_{\a\ob\oc e} \, \w^e_0 & = & 
    \d^\z_\b \, \left( \d^\r_\a \, \w^0_\oc - \d_{\a\c} \, \w^\r_{2m-1} 
      - \w^\r_{\a\c} \right) \, + \, 
    \d^\z_\c \, \left( \d^\r_\a \, \w^0_\ob - \d_{\a\b} \, \w^\r_{2m-1} 
      - \w^\r_{\a\b} \right)  \, .
  \label{eqn:sl8} \\
  \nonumber
\end{eqnarray}

The first three equations (\ref{eqn:sl1}, \ref{eqn:sl2}, \ref{eqn:sl3}) force 
$0 = \w^{\oz}_{m0} = \w^{\r}_{m0}$.  This is seen as follows.  Equations (\ref{eqn:sl2}, \ref{eqn:sl3}) imply that there exist functions $r^\r_e$, $r^\oz_e$ such that 
\begin{equation}
\label{eqn:sl9}
  r^{\r\z}_{\a(2m-1)^2e} = \d^\r_\a r^\oz_e
  \quad \hbox{ and } \quad
  r^{\r\z}_{\oa(2m-1)^2e} = \d^\z_\a r^\r_e \, .
\end{equation}
Immediately, the vanishing $0 = r^{\r\z}_{(2m-1)^{3}e}$ of (\ref{eqn:sl1}) 
yields $0 = r^{\rho}_{2m-1} = r^{\oz}_{2m-1}$.  Next, the symmetries 
$r^{\r\z}_{\a(2m-1)^2\b} = r^{\r\z}_{\b(2m-1)^2\a}$ and 
$r^{\r\z}_{\oa(2m-1)^2\ob} = r^{\r\z}_{\ob(2m-1)^2\oa}$ 
imply $\d^\r_\a r^\oz_\b = \d^\r_\b r^\oz_\a$ and 
$\d^\z_\a r^\r_\ob = \d^\z_\b r^\r_\oa$.  If 
$m-1 \ge 2$, we may pick $\rho=\alpha\not=\beta$ to see that $0 = r^\oz_\b$.  
Similarly, $0 = r^\r_\oa$.  Likewise, working with the symmetry 
$r^{\r\z}_{\a(2m-1)^2\ob} = r^{\r\z}_{\ob(2m-1)^2\a}$, we may deduce 
$0 = r^\r_\a = r^\oz_\ob$.  

At this point we have shown that the functions $r^{\oz}_{e}$ and $r^{\r}_{e}$ are identically zero.  It follows from (\ref{eqn:sl2}, \ref{eqn:sl3}, \ref{eqn:sl9}) that
\begin{equation}
\label{eqn:sl10}
  r^{\r\z}_{ab(2m-1)^2} = 0 \, , \quad \hbox{ and } \quad
  0 \ = \ \w^\r_{m0} \ = \ \w^\oz_{m0}  \, .
\end{equation}
\medskip
{\it Remark.}  The case $m-1 = 1$ is addressed separately in 
\S\ref{sec:SLsmalln}.
\medskip

\subsection{$F_{4,Y}$ -- the conormal directions $u = \r0, m\r$}
\label{sec:SLF4rmdr}

Equations (\ref{eqn:F4}, \ref{eqn:SLorder3}, \ref{eqn:SL4n0}, \ref{eqn:sl10}) 
yield
\begin{eqnarray}
  0 & = & 
  r^{\r0}_{(2m-1)^{3}e} \, \w^{e}_{0} \ = \ 
  r^{\r0}_{\oa\ob ce} \, \w^e_0 \ = \ 
  r^{\r0}_{\a\b\c e} \, \w^e_0 
  \label{eqn:sl11}  \\
  r^{\r0}_{\a(2m-1)^2e} \, \w^e_0 & = & -\half \, 
    \left( \w^\r_{\a0} + \d^\r_\a \, \w^0_{2m-1} \right)  
  \label{eqn:sl12}  \\
  r^{\r0}_{\oa(2m-1)^2e} \, \w^e_0 & = & \half \, \w^\r_{m\a}  
  \label{eqn:sl13}  \\
  r^{\r0}_{\a\b(2m-1)e} \, \w^e_0 & = & \half
    \big( \d^\r_\a \, \d_{\b\e} + \d^\r_\b \, \d_{\a\e} \big) \, 
    \left( \w^\oe_{2m-1} - \half \, \w^0_\e \right) 
  \label{eqn:sl14}  \\
  r^{\r0}_{\a\ob(2m-1)e} \, \w^e_0 & = & \half \, 
    \left( \w^\r_{\a\b} + 
    \left( \d^\r_\a \, \d_{\b\e} 
    + \d^\r_\e \, \d_{\a\b} \right) \, \w^\e_{2m-1}
    - \half \, \d^\r_\a \, \w^0_\ob \right)  
  \label{eqn:sl15}  \\
  r^{\r0}_{\a\b\oc e} \, \w^e_0 & = & 
    \half \, \d^\r_\a \, \left( \w^\ob_\oc + \w^\c_\b - \d_{\b\c} \, 
      \left( \w^0_0 + \w^{2m-1}_{2m-1} \right) \right)   
  \label{eqn:sl16}  \\
    & & + \ 
    \half \, \d^\r_\b \, \Big( \w^\oa_\oc + \w^\c_\a - \d_{\a\c} \, 
      \left( \w^0_0 + \w^{2m-1}_{2m-1} \right) \Big) \, ,   
      \nonumber 
\end{eqnarray}
and 
\begin{eqnarray}
  0 & = & r^{m\r}_{(2m-1)^3e} \, \w^e_0 \ = \ 
  r^{m\r}_{\a\b ce} \, \w^e_0 \ = \ 
  r^{m\r}_{\oa\ob\oc e} \, \w^e_0
  \label{eqn:sl17}  \\
  r^{m\r}_{\a(2m-1)^2e} \, \w^e_0 & = & \half \, \w^\ovr_{\a0}  
  \label{eqn:sl18}  \\
  r^{m\r}_{\oa(2m-1)^2e} \, \w^e_0 & = & \half \, 
    \left( \d^\r_\a \, \w^0_{2m-1} - \w^\ovr_{m\a} \right)  
  \label{eqn:sl19}  \\
  r^{m\r}_{\a\ob(2m-1)e} \, \w^e_0 & = & \half \, 
    \left( - \w^\ovr_{\a\b} + 
    \big( \d^\r_\b \, \d_{\a\e} 
    + \d^\r_\e \, \d_{\a\b} \big) \, \w^\oe_{2m-1}
    + \half \, \d^\r_\b \, \w^0_\a \right)  
  \label{eqn:sl20}  \\
  r^{m\r}_{\oa\ob(2m-1)e} \, \w^e_0 & = & \half
    \big( \d^\r_\a \, \d_{\b\e} + \d^\r_\b \, \d_{\a\e} \big) \, 
    \left( \w^\e_{2m-1} + \half \, \w^0_\oe \right) 
  \label{eqn:sl21}  \\
  r^{m\r}_{\a\ob\oc e} \, \w^e_0 & = & 
    \half \, \d^\r_\b \, \Big( \w^\oa_\oc + \w^\c_\a - \d_{\a\c} \, 
      \left( \w^0_0 + \w^{2m-1}_{2m-1} \right) \Big)   
  \label{eqn:sl22}  \\
    & & + \  
    \half \, \d^\r_\c \, \Big( \w^\oa_\ob + \w^\b_\a - \d_{\a\b} \, 
      \left( \w^0_0 + \w^{2m-1}_{2m-1} \right) \Big) \, .
      \nonumber
\end{eqnarray}

Let's consider the various expressions for 
$\d^{\r}_{\a} \, \d^{\z}_{\b} \, \w^{0}_{2m-1} + \half \, \d^{\z}_{\b} \, \w^{\r}_{\a0} 
  - \half \, \d^{\r}_{\a} \, \w^{\oz}_{m\b}$ given above.  From (\ref{eqn:sl5}), and then 
  (\ref{eqn:sl10}):
\begin{displaymath}
   \d^\r_\a \, \d^\z_\b \, \w^0_{2m-1} + \half \, \d^\z_\b \, \w^\r_{\a0} - \half \, \d^\r_\a \, \w^\oz_{m\b} 
   \ = \ r^{\r\z}_{\a\ob(2m-1)e} \, \w^e_{0}
   \ = \ r^{\r\z}_{\a\ob(2m-1)\e} \, \w^\e_{0} + r^{\r\z}_{\a\ob(2m-1)\oe} \, \w^\oe_{0} \, .
\end{displaymath}
And with (\ref{eqn:sl12}, \ref{eqn:sl19}) and (\ref{eqn:sl11}, \ref{eqn:sl17}):
\begin{displaymath}
\renewcommand{\arraystretch}{1.5}
\begin{array}{l}
  \d^\r_\a \, \d^\z_\b \, \w^0_{2m-1} + \half \, \d^\z_\b \, \w^\r_{\a0} - \half \, \d^\r_\a \, \w^\oz_{m\b}
  \ = \
  \left( 
    \d^{\r}_{\a} \, r^{m\z}_{\ob(2m-1)^{2}e} - \d^{\z}_{\b} \, r^{\r0}_{\a(2m-1)^{2}e}
  \right)  \, \w^{e}_{0}  
  \\  
  \hbox{\hspace{75pt}} = \  
  \left( 
    \d^{\r}_{\a} \, r^{m\z}_{\ob(2m-1)^{2}\e} - \d^{\z}_{\b} \, r^{\r0}_{\a(2m-1)^{2}\e}
  \right)  \, \w^{\e}_{0}
  \, + \,  
  \left( 
    \d^{\r}_{\a} \, r^{m\z}_{\ob(2m-1)^{2}\oe} - \d^{\z}_{\b} \, r^{\r0}_{\a(2m-1)^{2}\oe}
  \right)  \, \w^{\oe}_{0} \, .
\end{array}
\end{displaymath}
A comparison of these expressions yields 
\begin{eqnarray*}
  r^{\r\z}_{\a\ob(2m-1)\e} 
  & = & 
  \d^{\r}_{\a} \, r^{m\z}_{\ob(2m-1)^{2}\e} - \d^{\z}_{\b} \, r^{\r0}_{\a(2m-1)^{2}\e}
  \\
  r^{\r\z}_{\a\ob(2m-1)\oe}
  & = & 
  \d^{\r}_{\a} \, r^{m\z}_{\ob(2m-1)^{2}\oe} - \d^{\z}_{\b} \, r^{\r0}_{\a(2m-1)^{2}\oe} \, .
\end{eqnarray*}
The symmetry in $(\a,\e)$ on the left side of the first equation, and the symmetry in $(\ob,\oe)$ on the left side of the second equation force
\begin{equation}
\label{eqn:sl23}
  r^{m\z}_{\ob(2m-1)^{2}\e}  \  =  \  0  \  =  \  r^{\r0}_{\a(2m-1)^{2}\oe}  \, ,
\end{equation}
respectively.  This updates the formulas above to 
\begin{equation}
\label{eqn:sl24}
\renewcommand{\arraystretch}{1.5}
\begin{array}{rcrcl}
  \d^\r_\a \, \d^\z_\b \, \w^0_{2m-1} + \half \, \d^\z_\b \, \w^\r_{\a0} - \half \, \d^\r_\a \, \w^\oz_{m\b}
  & = & 
  r^{\r\z}_{\a\ob(2m-1)\e} \, \w^\e_{0} & + & r^{\r\z}_{\a\ob(2m-1)\oe} \, \w^\oe_{0} \\
  & = & 
  - \d^{\z}_{\b} \, r^{\r0}_{\a(2m-1)^{2}\e}  \, \w^{\e}_{0}
  & + & \d^{\r}_{\a} \, r^{m\z}_{\ob(2m-1)^{2}\oe}  \, \w^{\oe}_{0} \, .
\end{array}
\end{equation}
At this point we have 
\begin{equation}
\label{eqn:sl25}
  \half \, \w^{\r}_{m\a} \ = \ 
  r^{\r0}_{\oa(2m-1)^{2}e} \, \w^{e}_{0} \ = \ 
  0 \, .
\end{equation}
The first equality is just (\ref{eqn:sl13}).  The second equality is a consequence of 
(\ref{eqn:sl11}, \ref{eqn:sl23}).  Similarly, (\ref{eqn:sl17}, \ref{eqn:sl18}, \ref{eqn:sl23}) yield
\begin{equation}
\label{eqn:sl26}
  \half \, \w^{\ovr}_{\a0} \ = \
  r^{m\r}_{\a(2m-1)^{2}e} \, \w^{e}_{0} \ = \ 
  0 \, .
\end{equation}
These two equations, in conjunction with (\ref{eqn:sl4}, \ref{eqn:sl6}), yield
\begin{displaymath}
  r^{\r\z}_{\oa\ob(2m-1)e} \ = \ 0 \ = \ r^{\r\z}_{\a\b(2m-1)e} \, .
\end{displaymath}
This, and (\ref{eqn:sl10}), implies
\begin{displaymath}
  r^{\r\z}_{\a\ob(2m-1)e} = 0 \, .
\end{displaymath}
Now, from (\ref{eqn:sl24}), we may conclude
\begin{displaymath}
  r^{\r0}_{\a(2m-1)^{2}\e} \ = \ 0 \ = \ r^{m\z}_{\ob(2m-1)^{2}\oe} \, .
\end{displaymath}
The first equality and (\ref{eqn:sl11}, \ref{eqn:sl12}, \ref{eqn:sl23}) give us
\begin{equation}
\label{eqn:sl27}
  \w^{\r}_{\a0} + \d^{\r}_{\a} \, \w^{0}_{2m-1} \ = \ 0 \, .
\end{equation}
Similarly, the second equality and (\ref{eqn:sl17}, \ref{eqn:sl19}, \ref{eqn:sl23}) yield
\begin{equation}
\label{eqn:sl28}
  \w^{\ovr}_{m\a} - \d^{\r}_{a} \, \w^{0}_{2m-1} \ = \ 0 \, .
\end{equation}
Finally, with (\ref{eqn:sl5}, \ref{eqn:sl27}, \ref{eqn:sl28}) we have $r^{\r\z}_{\a\ob(2m-1)e} = 0$.

Let's pause for a moment to assess our progress toward showing that the coefficients 
$r^{u}_{abce}$ vanish ($u = \r\z, \r0, n\r$).
\begin{itemize}
  \item
  The vanishing of the coefficients in (\ref{eqn:sl2}--\ref{eqn:sl6}) is now ensured by 
  (\ref{eqn:sl10}, \ref{eqn:sl25}--\ref{eqn:sl28}).  Notice in particular that the only potentially non-zero coefficients corresponding to the conormal direction $u = \r\z$ are $r^{\r\z}_{\a\b\oc\oe}$.
  \item
  The vanishing of the coefficients in (\ref{eqn:sl12}, \ref{eqn:sl13}) is given by 
  (\ref{eqn:sl25}, \ref{eqn:sl27}).
  \item  
  The vanishing of the coefficients in (\ref{eqn:sl18}, \ref{eqn:sl19}) is equivalent to 
  (\ref{eqn:sl26}, \ref{eqn:sl28}).
\end{itemize}
It remains to address the coefficients appearing in 
(\ref{eqn:sl7}, \ref{eqn:sl8}, \ref{eqn:sl14}--\ref{eqn:sl16}, \ref{eqn:sl20}--\ref{eqn:sl22}).

The next portion of the analysis focuses on the equations (\ref{eqn:sl16}, \ref{eqn:sl22}).
They both (individually) imply that there are functions $r_{\b\oc e}$ such that 
\begin{displaymath}
  \w^{\ob}_{\oc} + \w^{\c}_{\b} - \d_{\b\c} \, 
  \left(  
    \w^{0}_{0} + \w^{2m-1}_{2m-1}
  \right)
  \ = \ 
  r_{\b\oc e} \, \w^{e}_{0} \, .
\end{displaymath}
In particular, 
\begin{displaymath}
  r^{\r0}_{\a\b\oc e} = \half \, \d^{\r}_{\a} \, r_{\b\oc e} + \half \, \d^{\r}_{\b} \, r_{\a\oc e}
  \quad \hbox{ and } \quad
  r^{m\r}_{\a\ob\oc e} = \half \, \d^{\r}_{\b} \, r_{\a\oc e} + \half \, \d^{\r}_{\c} \, r_{\a\ob e} \, .
\end{displaymath}
The symmetry of $r^{\r0}_{\a\b\oc\oe}$ in $(\oc, \oe)$,  and the symmetry of 
$r^{m\r}_{\a\ob\oc\e}$ in $(\a, \e)$ imply $r_{\b\oc\oe}$ is symmetric in $(\oc, \oe)$, and 
$r_{\a\oc\e}$ is symmetric in $(\a, \e)$, respectively.

Now the symmetry of $r^{\r0}_{\a\b\oc\e}$ in $(\b,\e)$ yields $r_{\a\oc\e} = 0$; and the symmetry of $r^{m\r}_{\a\ob\oc\oe}$ in $(\oc, \oe)$ yields $r_{\a\ob\oc} = 0$.  It is now a consequence of the equations (\ref{eqn:sl11}, \ref{eqn:sl16}, \ref{eqn:sl17}, \ref{eqn:sl22}) that 
\begin{displaymath}
  r^{\r0}_{abce} \, , \ r^{m\r}_{abce} \ = \ 0
  \quad \hbox{ if } a, b, c, e \not= 2m-1 \, .
\end{displaymath}
In fact, the only remaining, potentially non-zero, coefficients $r^{\r0}_{abce}$, 
$r^{m\r}_{abce}$ are $r^{\r0}_{\a\b\oc(2m-1)}$ and $r^{m\r}_{\a\ob\oc(2m-1)}$.

This brings us to the final stage of our analysis of the coefficients of $F_{4,Y}$.  From $(\ref{eqn:sl14})$, we see that there are functions $r_{\a\oe}$ so that 
\begin{displaymath}
  \w^\oa_{2m-1} - \half \, \w^0_\a = r_{\a\oe} \, \w^\oe_0 \, .
\end{displaymath}
In particular, (\ref{eqn:sl20}) implies
\begin{displaymath}
  \left(
    r^{m\r}_{\a\ob(2m-1)\oe} - \half \, \d^\r_\b \, r_{\a\oe} 
  \right) \, \w^\oe_0
  \ = \ 
  \half \, \left(
    -\w^\ovr_{\a\b} + \d_{\a\b} \, \w^\ovr_{2m-1} + \d^\r_\b \, \w^0_\a
  \right) \, ;
\end{displaymath}
which, with (\ref{eqn:sl7}), allows us to write
\begin{displaymath}
  r^{\r\z}_{\a\b\oc\oe} \ = \ 
  2 \, \d^\r_\a \, \left(
    r^{m\z}_{\b\oc(2m-1)\oe} - \half \, \d^\z_\c \, r_{\b\oe}
  \right)
  \ + \ 2 \, \d^\r_\b \, \left(
    r^{m\z}_{\a\oc(2m-1)\oe} - \half \, \d^\z_\c \, r_{\a\oe}
  \right) \, .
\end{displaymath}
The symmetry in $(\oc, \oe)$ on the right forces $r_{\a\oe} = 0$, and we may 
conclude 
\begin{eqnarray}
  \label{eqn:sl30}
  0 & = & 
  \w^\oa_{2m-1} - \half \, \w^0_\a 
  \\
  \label{eqn:sl30a}
  \half \, r^{\r\z}_{\a\b\oc\oe} & = & 
  \d^\r_\a \, r^{m\z}_{\b\oc(2m-1)\oe} 
  \ + \ \d^\r_\b \, r^{m\z}_{\a\oc(2m-1)\oe} \, .
\end{eqnarray}
The analogous argument with (\ref{eqn:sl8}, \ref{eqn:sl15}, \ref{eqn:sl21})
yields 
\begin{equation}
\label{eqn:sl31}
  0 \ = \ \w^\a_{2m-1} + \half \, \w^0_\oa \, .
\end{equation}

With (\ref{eqn:sl14}, \ref{eqn:sl30}) we deduce that 
$r^{\r0}_{\a\b(2m-1)e} = 0$, and we may now conclude that all the coefficients 
$r^{\r0}_{abce}$ vanish.  In particular, (\ref{eqn:sl16}) yields
\begin{equation}
\label{eqn:sl29}
  \w^\a_\b + \w^\ob_\oa = \d_{\a\b} \, 
  \left( \w^0_0 + \w^{2m-1}_{2m-1} \right) \, .
\end{equation}

Similarly, the vanishing of $r^{m\r}_{abce}$ is a 
consequence of (\ref{eqn:sl21}, \ref{eqn:sl31}).  And this, along with 
(\ref{eqn:sl30a}), yields $r^{\r\z}_{abce} = 0$.  Therefore the only nonzero 
coefficients of $F_{4,Y}$ are given by (\ref{eqn:SLF4}) and we have established
\begin{displaymath}
  F_{4,Y} = F_{4,X} \, .
\end{displaymath}
Finally, note that $r^{\r\z}_{abce}=0$ and 
(\ref{eqn:sl7}, \ref{eqn:sl8}, \ref{eqn:sl30}, \ref{eqn:sl31})
provide us with 
\begin{equation}
\label{eqn:sl32}
\renewcommand{\arraystretch}{1.5}
\begin{array}{rcl}
  \w^\r_{\a\b} & = & \d^\r_\a \, \w^0_\ob + \half \, \d_{\a\b} \, \w^0_\ovr \\
  \w^\oz_{\a\b} & = & \d^\z_\b \, \w^0_\a + \half \, \d_{\a\b} \, \w^0_\z \, . 
\end{array}
\end{equation}

\subsection{Computations of $F_{k,Y}$, $k\ge5$}
\label{sec:SLF5+}

In this section I will show that the higher order invariants $F_{k,Y}$, 
$k \ge 5$, vanish on 
$\F^4_U$ (the sub-bundle of $\F^3_U$ on which the normalizations of \ref{sec:SLF4n0} hold).  We begin with the coefficients of $F_{5,Y}$, which are given by Proposition \ref{prop:F}.  I will use the relations (\ref{eqn:SLorder3}, \ref{eqn:SL4n0}, \ref{eqn:sl10}, \ref{eqn:sl25}--\ref{eqn:sl32}) without mention.

Start with the coefficients corresponding to the conormal direction 
$u = \r\z$; we will see that the $\w^0_v$ vanish.  First, 
\begin{eqnarray*}
  0  & = &  
  r^{\r\z}_{(2m-1)^3ae} \, \w^e_0 \ = \ 
  r^{\r\z}_{(2m-1)^2\a\b e} \, \w^e_0  \ = \ 
  r^{\r\z}_{(2m-1)^2\oa\ob e} \, \w^e_0 \\
  \half \, \d^\r_\a \, \d^\z_\b \, \w^0_{m0} 
  & = &
  r^{\r\z}_{(2m-1)^2\a\ob e} \, \w^e_0 \, . 
\end{eqnarray*}
Symmetry in the lower indices of $r^u_{abcde}$ (and $m-1 > 1$) forces
$$ \w^0_{m0} = 0 \, .$$

Next, 
\begin{eqnarray*}
  r^{\r\z}_{(2m-1)\a\b\c e} \, \w^e_0  \ = \
  r^{\r\z}_{(2m-1)\oa\ob\oc e} \, \w^e_0 
  & = &  0  \\
  r^{\r\z}_{(2m-1)\a\b\oc e} \, w^e_0 
  & = & 
  \half \, \d^\z_\c \, \left(
    \d^\r_\a \, \w^0_{\b0} + \d^\r_\b \, \w^0_{\a0} 
  \right)  \\
  r^{\r\z}_{(2m-1)\a\ob\oc e} \, w^e_0 
  & = & 
  - \half \, \d^\r_\a \, \left(
    \d^\z_\b \, \w^0_{m\c} + \d^\z_\c \, \w^0_{m\b} 
  \right)  \, .
\end{eqnarray*}
As before, symmetry forces
$$ 0 \ = \ \w^0_{\a0} \ = \ \w^0_{m\b} \, .$$

Finally, 
\begin{eqnarray*}
  0  &  =  &
  r^{\r\z}_{\a\b\c\vartheta e} \, \w^e_0  \ = \ 
  r^{\r\z}_{\a\b\c\overline{\vartheta} e} \, \w^e_0  \ = \ 
  r^{\r\z}_{\a\ob\oc\overline{\vartheta} e} \, \w^e_0  \ = \ 
  r^{\r\z}_{\oa\ob\oc\overline{\vartheta} e} \, \w^e_0  \\
  r^{\r\z}_{\a\b\oc\overline{\vartheta} e} \, \w^e_0
  &  =  &
  - \d^\r_\a \, \d^\z_\c \, \w^0_{\b\vartheta}
  - \d^\r_\a \, \d^\z_\vartheta \, \w^0_{\b\c}
  - \d^\r_\b \, \d^\z_\c \, \w^0_{\a\vartheta}
  - \d^\r_\b \, \d^\z_\vartheta \, \w^0_{\a\c} \, ,
\end{eqnarray*}
and once again symmetry forces
$$ \w^0_{\a\b} \ = \ 0 \, .$$
We conclude that the coefficients of $F_{5,Y}$ corresponding to the conormal 
direction $u = \r\z$ vanish.  

Straightforward, if lengthly, computations show 
that the remaining coefficients vanish as well.  ({\it Remark}.  Neither these computations, nor those that follow, require $m - 1 > 1$.)

Additional calculations with Proposition \ref{prop:F} yield 
\begin{center}
  $F_{6,Y} \, , \ F_{7,Y} \, , \ F_{8,Y} \, , \ F_{9,Y} \ = \ 0 \, ,$
\end{center}
completing the proof of Theorem \ref{thm:SL} (in the case $m > 2$).

\subsection{When \boldmath$m=2$\unboldmath.}
\label{sec:SLsmalln}

Since $m=2$, we have $\a = 1$ and $\oa = \overline 1 = 2$.  For consistency I 
will continue to use the notation $\a, \oa$, rather than $1, 2$, but will 
abbreviate $2m-1 = 3$.  As before $1 \le c,e \le 2m-1 = 3$.  To complete the 
proof of Theorem \ref{thm:SL} we need to do two things:  
\begin{itemize}
  \item[(1)] 
    show that $r^{u}_{abce} = 0$ for $u = \a0, m\a, \a\a$, 
    and 
  \item[(2)] 
    show that $\w^{0}_{u} = 0$ for $u = \a0, m\a, \a\a, m0$.
\end{itemize}
First, (\ref{eqn:F4}) yields
\begin{displaymath}
  0 \ = \  r^{\a0}_{\a\a\a e} \ = \ r^{\a0}_{\oa\oa ce}  \ = \ 
    r^{m\a}_{\a\a ce} \ = \ r^{m\a}_{\oa\oa\oa e}  \ = \ 
    r^{\a\a}_{\a\a\a e} \ = \ r^{\a\a}_{\oa\oa\oa e} \ = \ r^{\a\a}_{333e} \, ,
\end{displaymath}
\begin{displaymath}
\renewcommand{\arraystretch}{1.5}
\begin{array}{r@{ \ \ = \ \ }r@{ \ \ = \ \ }l}
  \w^{\a}_{m0} 
  & 2 \, r^{\a\a}_{\oa 33e} \, \w^{e}_{0} 
  & -\frac{4}{3} \, r^{\a0}_{333e} \, \w^{e}_{0} \\
  \w^{\oa}_{m0}
  & 2 \, r^{\a\a}_{\a33e} \, \w^{e}_{0}
  & \frac{4}{3} \, r^{m\a}_{333e} \, \w^{e}_{0} \\
  \w^{\oa}_{\a0}  &
  r^{\a\a}_{\a\a3e} \, \w^{e}_{0}  & 
  2 \, r^{m\a}_{\a33e} \, \w^{e}_{0}  \\
  \w^{0}_{3} + \half \, \w^{\a}_{\a0} - \half \, \w^{\oa}_{m\a}   & 
  r^{\a\a}_{\a\oa3e} \, \w^{e}_{0}  & 
  r^{m\a}_{\oa33e} \, \w^{e}_{0} - r^{\a0}_{\a33e} \, \w^{e}_{0}  \\
  \w^{\a}_{m\a}  & 
  - r^{\a\a}_{\oa\oa3e} \, \w^{e}_{0}  &
  2 \, r^{\a0}_{\oa33e} \, \w^{e}_{0}  \\
  \w^0_\a + \w^\oa_3 - \w^\oa_{\a\a}  & 
  \half \, r^{\a\a}_{\a\a\oa e} \, \w^e_0  &
  2 \, r^{m\a}_{\a\oa3e} \, \w^e_0 - r^{\a0}_{\a\a3e} \, \w^e_0  \\
  \w^0_\oa - \w^\a_3 - \w^\a_{\a\a}  & 
  \half \, r^{\a\a}_{\a\oa\oa e} \, \w^e_0  & 
  -2 \, r^{\a0}_{\a\oa3e} \, \w^e_0 + r^{m\a}_{\oa\oa3e} \, \w^e_0  \\ 
  \w^\a_\a + \w^\oa_\oa - \w^0_0 - \w^3_3  &  
  r^{m\a}_{\a\oa\oa e} \, \w^e_0  & 
  r^{\a0}_{\a\a\oa e} \, \w^e_0  \, .
\end{array}
\end{displaymath}
These relations force the coefficients to vanish, establishing (1). 

Next calculations with Proposition \ref{prop:F} produce
\begin{displaymath}
  0 \ = \ r^{\a\a}_{333ce} \, \w^e_0 \ = \ 
  r^{\a\a}_{\a\a33e} \, \w^e_0 \ = \ 
  r^{\a\a}_{\oa\oa33e} \, \w^e_0 \, , \quad \hbox{ and } \quad 
  r^{\a\a}_{\a\oa33e} \, \w^e_0 \ = \ \fourth \w^0_{m0} \, .
\end{displaymath}
Symmetry in the lower indices of $r^{\a\a}_{abcde}$ forces 
\begin{displaymath}
  \w^0_{m0} \ = \ 0 \, .
\end{displaymath}
Additional computations yield
\begin{eqnarray*}
  \w^0_{\a0} & = & 
  -r^{\a0}_{\a\a33\oa} \, \w^\oa_0 \ = \ 
  2 \, r^{m\a}_{\a\oa33\oa} \, \w^\oa_0 \ = \ 
  r^{\a\a}_{\a\a\oa3\oa} \, \w^\oa_0 \ = \ 
  - \textstyle{\frac43} \, r^{m0}_{\a333\oa} \, \w^\oa_0 \\
  \w^0_{m\a} & = & 
  2 \, r^{\a0}_{\a\oa33\a} \, \w^\a_0 \ = \ 
  - r^{m\a}_{\oa\oa33\a} \, \w^\a_0 \ = \ 
  - r^{\a\a}_{\a\oa\oa3\a} \, \w^\a_0 \ = \ 
  \textstyle{\frac43} \, r^{m0}_{\oa333\a} \, \w^\a_0 \\
  \w^0_{\a\a} & = & 
  2 \, r^{\a0}_{\a\a\oa33} \, \w^3_0 \ = \ 
  - r^{m\a}_{\a\oa\oa33} \, \w^3_0 \ = \ 
  -\fourth \, r^{\a\a}_{\a\a\oa\oa3} \, \w^3_0 \ = \ 
  2 \, r^{m0}_{\a\oa333} \, \w^3_0 \, .
\end{eqnarray*}
And again the coefficients, and therefore the $\w^0_u$, must vanish.

\section{Concluding remarks}
\label{sec:conclusions}

The Veronese embedding $v_2(\P^{2m-1})$ of $\P^{2m-1}$ may be identified with 
the adjoint variety of $\mathrm{Sp}_{2m}\C$, which is known to be rigid to 
order three (c.f. \S\ref{sec:veronese}).  So it is natural to ask if 
the adjoint varieties of the simple Lie groups are all rigid to order three.
\vspace{0.1in}\\
There is some reason to hope that this is the case, as there are many similarities amongst these spaces:  Given the adjoint variety
of a simple Lie group it is the case that 
\begin{itemize}
  \item The $F_k$ may be normalized to zero, $k \ge 5$.
  \item On the reduced frame bundle there is a single Fubini quartic 
        $\P |F_4| \in S^4 \, T_{1,x}^*$.
        Here, as in \S\ref{sec:SLFk}, $T_{1,x} \subset T_xX$ is a contact 
        hyperplane.
  \item The Fubini cubics $|F_3| \subset S^3 \, T_{1,x}$ are the derivatives 
        of $|F_4|$.
  \item $\mathcal C_{2,x} = \mathcal C_x \subset \P T_{1,x}$. 
\end{itemize}
The adjoint varieties of $\mathrm{SL}_{m+1}\C$ and $\mathrm{Sp}_{2m}\C$ are 
degenerate in the following sense.  For $\mathrm{Sp}_{2m}\C$, that single 
Fubini quartic is zero, and $\mathcal C_{2,x} = \emptyset$ 
(cf. \S\ref{sec:veronese} 
and \cite{L99}).  In the case of of $\mathrm{SL}_{m+1}\C$, 
$\mathcal C_{2,x} = \P^{m-1} \sqcup \P^{m-1} \subset \P \, T_{1,x}$ is the 
disjoint union of two linear spaces, and the Fubini quartic factors as the 
square of two quadrics.  (See \S\ref{sec:SLFk}.)  The adjoint representation 
fails to be fundamental for these groups.  

The adjoint representation is fundamental for the remaining simple Lie groups.
As a consequence, 
\begin{itemize}
  \item $\mathcal C_{2,x}$ is a generalized minuscule variety, and the closed 
    orbit of a semi-simple $H\subset G$ in $\P \, T_{1,x}$.
  \item The Fubini quartic is irreducible, and its zero locus in $\P\, T_{1,x}$
    is the tangential variety of $\mathcal C_{2,x}$.
\end{itemize}
See \cite{LM02} for details.



\begin{thebibliography}{AA99}

\bibitem{BCG3} R. Bryant, S.-S. Chern, R.B. Gardner, H. Goldschmidt \&
P. Griffiths, {\it Exterior Differential Systems}, MSRI Publications, 
Springer, 1990.

\bibitem{FH91}
W. Fulton \& J. Harris, {\it Representation Theory: A First Course}, Springer-Verlag, 
New York, 1991.

\bibitem{F18} G. Fubini,   {\it Studi relativi all'elemento
lineaaare proiettivo di una ipersuficie}, Rend. Acad. Naz. dei
Lincei (1918) 99--106.

\bibitem{GH79}  P. Griffiths \& J. Harris, {\it Algebraic geometry and local 
differential geometry},  Ann. Sci. Ecole Norm. Sup. {\bf 12} (1979) 355--432.

\bibitem{HM06}
J.-M. Hwang \& N. Mok,
{\it Rigidity of irreducible Hermitian symmetric spaces of the compact type 
under K\"ahler deformation}, Invent. Math. {\bf 131} (1998), no. 2, 393--418.

\bibitem{HY03}
J.-M. Hwang \& L. Yamaguchi,  
{\it Characterization of Hermitian symmetric spaces by fundamental forms}.  Duke Math. J.  
{\bf 120}  (2003),  no. 3, 621--634. 

\bibitem{IL03} T.A. Ivey \& J.M. Landsberg, {\it Cartan for beginners: differential
geometry via moving frames and exterior differential systems}, 
Graduate Studies in Mathematics 61, A.M.S. (2003).

\bibitem{K61}
B. Kostant, 
{\it Lie algebra cohomology and the generalized Borel--Weil theorem},
Ann. of Math. (2) {\bf 74} (1961) 329--387, MR0142696.

\bibitem{L96} J.M. Landsberg, {\it Differential-geometric characterizations of 
complete intersections}, J. Diff. Geom. {\bf 44} (1996), 32-73.

\bibitem{L99}
J.M. Landsberg, 
{\it On the infinitesimal rigidity of homogeneous varieties}, Compositio Math. 
{\bf 118} (1999) 189--201.

\bibitem{L06}
J.M. Landsberg,
{\it Griffiths-Harris rigidity of compact Hermitian symmetric spaces},
to appear in J. Diff. Geom.

\bibitem{LM02}
J.M. Landsberg \& L. Manivel,
{\it Classification of simple Lie algebras via projective geometry},
Selecta Math. {\bf 8} (2002) 137--159, MR1890196.

\bibitem{LM03}
J.M. Landsberg \& L. Manivel, 
{\it On the projective geometry of rational homogeneous varieties},
Comment. Math. Helv. {\bf 78} (1) (2003) 65--100, MR1966752.

\bibitem{S88}
Y. Se-ashi, 
{\it On differential invariants of integrable finite type linear differential equations},
Hokkaido Math. J. {\bf 17} (1988), no. 2, 151--195, MR0945853.

\end{thebibliography}
\end{document}